\input amstex
\input amsppt.sty
\magnification=1200 \baselineskip=24truept
 \NoBlackBoxes
 \hsize=16truecm
 \vsize=22.5truecm
\def\q{\quad}

\def\({\left(}
\def\){\right)}
\def\[{\left[}
\def\]{\right]}
\def\mod#1{\ (\text{\rm mod}\ #1)}
\def\t{\text}
\TagsOnRight
\def\f{\frac}

\def\e{\equiv}

\def\ord{\text{\rm ord}}
\def\trm#1{(\text{\rm #1})}

\def\inZ{\in\Bbb Z}
\def\p3{(p-(\f p3))/3}
\def\up3{u_{\f{p-(\f p3)}3}}
 \def\vp3{v_{\f{p-(\f p3)}3}}

\let \pro=\proclaim
\let \endpro=\endproclaim

\topmatter
\title Congruences concerning Lucas' law of repetition\endtitle
\author Zhi-Hong Sun\endauthor
\affil School of Mathematical Sciences, Huaiyin Normal University,
\\ Huaian, Jiangsu 223001, PR China
\\ E-mail: zhihongsun$\@$yahoo.com
\\ Homepage: http://www.hytc.edu.cn/xsjl/szh\endaffil
\nologo \NoRunningHeads \abstract{ Let $P,Q\in\Bbb Z$, $U_0=0,\
U_1=1$ and $U_{n+1}=PU_n-QU_{n+1}$. In this paper we obtain a
general congruence for $U_{kmn^r}/U_k\pmod {n^{r+1}}$, where
$k,m,n,r$ are positive integers. As applications we extend Lucas'
law of repetition and characterize the square prime factors of $
a^n+1$ or $S_n$, where $\{S_n\}$ is given by $S_1=P^2+2$ and
$S_{k+1}=S_k^2-2\ (k\ge 1)$.
\par\q
}
\endabstract
\endtopmatter
\document
\baselineskip=24truept
 \subheading{1. Introduction}
\par For complex numbers $P$ and $Q$ the Lucas sequences $\{U_n(P,Q)\}$ and $\{V_n(P,Q)\}$
are defined by
$$U_0(P,Q)=0,\ U_1(P,Q)=1,\ U_{n+1}(P,Q)=PU_n(P,Q)
-QU_{n-1}(P,Q)\ (n\ge 1)$$ and
$$V_0(P,Q)=2,\ V_1(P,Q)=P,\ V_{n+1}(P,Q)=PV_n(P,Q)-
QV_{n-1}(P,Q)\ (n\ge 1).$$
\par Let $k,r,s$ be positive integers, $P,Q\inZ$  and
$U_k=U_k(P,Q)$. If $p$ is a prime such that $p\nmid Q$, $p^s\not=2$,
$p^s\mid U_k$ and $p^{s+1}\nmid U_k$, then $p^{r+s}\mid U_{kp^r}$
but $p^{r+s+1}\nmid U_{kp^r}$. This result is called  Lucas' law of
repetition (see [1,3,4]) since Lucas first obtained it for the
Fibonacci sequence $F_n=U_n(1,-1)$.
\par In this paper we mainly obtain a
general congruence for $U_{kmn^r}/U_k\mod {n^{r+1}}$, where
$k,m,n,r$ are positive integers. Using this congruence we extend
Lucas' law of repetition and characterize the square prime factors
of $ a^n+1$ or $S_n$, where $\{S_n\}$ is given by $S_1=P^2+2$ and
$S_{k+1}=S_k^2-2\ (k\ge 1)$.
\par Throughout the paper we let $(m,n)$ be the greatest common
divisor of integers $m$ and $n$, and let $(\f dp)$ be the Legendre
symbol. For given prime $p$ and positive integer $n$ let $\ord_pn$
denote the unique nonnegative integer $t$ satisfying $p^t\mid n$ and
$p^{t+1}\nmid n$.
 \subheading{2. Basic lemmas for Lucas
sequences}\par Let $U_n=U_n(P,Q)$ and $V_n=V_n(P,Q)$ be the Lucas
sequences given in Section 1. It is well known that
$$U_n=\f 1{\sqrt{P^2-4Q}}\Big\{\Big(\f{P+\sqrt{P^2-4Q}}2\Big)^n-
\Big(\f{P-\sqrt{P^2-4Q}}2\Big)^n\Big\}(P^2-4Q\not=0)\tag 1$$ and
$$V_n=\Big(\f{P+\sqrt{P^2-4Q}}2\Big)^n+
\Big(\f{P-\sqrt{P^2-4Q}}2\Big)^n.\tag 2$$  From (1) and (2) one can
easily  check the following known facts (see [1,3,4]):
$$\align &U_{2n}=U_nV_n,\
V_{2n}=V_n^2-2Q^n,\tag 3
\\&V_n^2-(P^2-4Q)U_n^2=4Q^n,\tag 4
\\&U_{n+k}=V_kU_n-Q^kU_{n-k}\ (n\ge k).\tag 5\endalign$$
By (5), if $U_k\not=0$, then
$U_{k(n+1)}/U_k=V_kU_{kn}/U_k-Q^kU_{k(n-1)}/U_k.$ Thus
$$U_{kn}/U_k=U_n(V_k,Q^k).\tag 6$$
\pro{Lemma 1} Suppose $P,Q,n\inZ$, $n\ge 0$ and $P^2-4Q\not=0$.
Then
$$U_n(P,Q)\e \cases nQ^{\f{n-1}2}\mod {|P^2-4Q|}&\t{if $2\nmid
n$,}\\ \f n2PQ^{\f{n-2}2}\mod {|P^2-4Q|}&\t{if $2\mid n$.}
\endcases$$
\endpro
Proof. We prove the lemma by induction. Clearly the result holds
for $n=0,1$. Now assume $n\ge 1$ and the result is true for $m\le
n$. Set $U_k=U_k(P,Q)$. By the inductive hypothesis we have
$$\aligned U_{n+1}&=PU_n-QU_{n-1}\\&\e
\cases P\cdot
nQ^{\f{n-1}2}-Q\cdot\f{n-1}2PQ^{\f{n-1}2-1}=\f{n+1}2PQ^{\f{n-1}2}
\mod {|P^2-4Q|}&\t{if $2\nmid n$,}
\\P\cdot \f n2PQ^{\f{n-2}2}-Q\cdot (n-1)Q^{\f{n-2}2}\e (n+1)Q^{\f
n2} \mod {|P^2-4Q|}&\t{if $2\mid n$.}
 \endcases\endaligned$$
So the result holds for $m=n+1$. Hence the lemma is proved.
  \pro{Lemma 2} If $P,Q\in\Bbb Z$,
$n,k\ge 0$, $U_n=U_n(P,Q)$ and $U_{2k}\not=0$. Then
$$\f{U_{(2n+1)k}}{U_k}\e (2n+1)Q^{kn}\mod{|P^2-4Q|U_k^2}$$
and
$$\f{U_{2kn}}{U_{2k}}\e nQ^{k(n-1)}\mod{|P^2-4Q|U_k^2}.$$
\endpro
Proof. Let $V_m=V_m(P,Q)$. From (4),(6) and Lemma 1 we have
$$\f{U_{(2n+1)k}}{U_k}=U_{2n+1}(V_k,Q^k)\e (2n+1)Q^{kn}\mod{|P^2-4Q|U_k^2}$$
and
$$\f{U_{2kn}}{U_{2k}}=U_n(V_{2k},Q^{2k})\e
\cases
 nQ^{2k\cdot \f{n-1}2}\mod{|P^2-4Q|U_{2k}^2}&\t{if $2\nmid n$,}
 \\ \f n2Q^{2k(\f n2-1)}V_{2k}\mod{|P^2-4Q|U_{2k}^2}&\t{if $2\mid
 n$.}\endcases $$
Note that $U_{2k}=U_kV_k$ and
$V_{2k}=V_k^2-2Q^k=(P^2-4Q)U_k^2+2Q^k$ by (3) and (4). By the
above we obtain the result.
 \pro{Lemma 3} Suppose $P,Q\in\Bbb Z$, $PQ(P^2-Q)(P^2-4Q)\not=0$
and $(P,Q)=1$. Then $U_k(P,Q)\not=0$ for any positive integer $k$.
\endpro
Proof. If $U_k(P,Q)=0$ for some positive integer $k$, applying (1)
we see that
$$\Big(\f{P^2-2Q+P\sqrt{P^2-4Q}}{2Q}\Big)^k=\Big(\f{P+\sqrt{P^2-4Q}}
{P-\sqrt{P^2-4Q}}\Big)^k=1.$$ If $P^2-4Q>0$, we must have
$P^2-2Q+P\sqrt{P^2-4Q}=\pm 2Q$. This yields $PQ=0$, which
contradicts the condition. If $P^2-4Q<0$, from the above we see
that
$$\Big(\f{P^2-2Q}{2Q}+\f{P\sqrt{4Q-P^2}}{2Q}i\Big)^k=1.$$
So there exists an integer $r\in\{0,1,\ldots,k-1\}$ such that
$$\f{P^2-2Q}{2Q}+\f{P\sqrt{4Q-P^2}}{2Q}i=\t{cos}\f{2\pi
r}k+i\t{sin}\f{2\pi r}k$$ and hence $\f{P^2-2Q}{2Q}=\t{cos}\f{2\pi
r}k$. But this is impossible since $Q\not= P^2$ and $(P,Q)=1$. Hence
the lemma is proved.
\pro{Lemma 4} Let $P,Q\in\Bbb Z$ and $(P,Q)=1$.
If $n$ is odd such that $n\mid U_k(P,Q)$ for some positive integer
$k$, then $(n,QV_k(P,Q))=1$.\endpro
Proof. Let $U_k=U_k(P,Q)$ and
$V_k=V_k(P,Q)$. From (4) we know that $V_k^2=(P^2-4Q)U_k^2+4Q^k$.
Since $n\mid U_k$ we get $V_k^2\e 4Q^k\mod n$. Suppose $p$ is a
prime divisor of $n$.
 If $p\mid Q$, then we must have
$p\mid V_k$ and $V_k\e P^k\mod p$ by (2). So $p\mid P$ and hence
$p\mid (P,Q)$. This contradicts the assumption $(P,Q)=1$. Hence
$p\nmid Q$ and so $p\nmid V_k$. Thus $(n,QV_k)=1$. \pro{Lemma 5 (see
[3,4])} Suppose $P,Q\in\Bbb Z$, $PQ(P^2-4Q)\not=0$, $(P,Q)=1$ and
$U_n=U_n(P,Q)$.
\par $\trm{i}$ If $p$ is an odd prime such that $p\nmid Q$, then $p\mid U_{p-(\f{P^2-4Q}p)}$.
\par $\trm{ii}$ If $r(p)$ is the least positive integer $n$ such that $p\mid U_n$, then
$p\mid U_m$ if and only if $r(p)\mid m$.
\endpro
\subheading{3. Main results}
 \pro{Theorem 1} Let $k,m,n,r$ be positive integers, $P,Q\in\Bbb Z$,
 $(P,Q)=1$, $PQ(P^2-Q)(P^2-4Q)\not=0$, $U_n=U_n(P,Q)$ and
$V_n=V_n(P,Q)$. If $n$ is a divisor of $U_k$, then
$$\f{U_{kmn^r}}{U_k}\e\cases
mn^rQ^{\f{k(mn^r-1)}2}\mod{n^{r+1}A}&\t{if $2\nmid
mn$}, \\
\f{mV_k}2n^rQ^{\f{k(mn^r-2)}2}\mod{n^{r+1}A}&\t{if $2\mid mn$},
\endcases $$ where
$A=|P^2-4Q|(\f{U_k}n)^2$.
\endpro
 Proof. From Lemma 3 we know that $U_t\not=0$
for $t\ge 1$. If $2\nmid n$, by Lemma 2 we have
$$\f{U_{kmn^{s+1}}}{U_{kmn^{s}}}\e nQ^{kmn^{s}\cdot\f{n-1}2}
\mod {|P^2-4Q|U_{kmn^s}^2}\q (s\ge 0).$$  Since $U_k\mid
U_{kmn^s}$ we find $U_{kmn^{s+1}}/U_{kmn^{s}}\e
nQ^{kmn^{s}\cdot\f{n-1}2} \mod {n^2A}$ for $s\ge 0$. For $1\le
j\le r$ it is clear that $n^{r+1}A\mid (n^2A)^jn^{r-j}$. So
$$\f{U_{kmn^r}}{U_{km}}=\prod_{s=0}^{r-1}\f{U_{kmn^{s+1}}}{U_{kmn^s}}
\e
\prod_{s=0}^{r-1}\big(nQ^{\f{km(n-1)n^s}2}\big)=n^rQ^{\f{km(n^r-1)}2}
\mod {n^{r+1}A}.$$
 Applying Lemma 2 and the above we then get
 $$\aligned\f{U_{kmn^r}}{U_k}&=\f{U_{km}}{U_k}\cdot\f{U_{kmn^r}}{U_{km}}
\\&\e\cases mQ^{\f{k(m-1)}2}\cdot
n^rQ^{\f{km(n^r-1)}2}=mn^rQ^{\f{k(mn^r-1)}2}
 \mod{n^{r+1}A}&\t{if $2\nmid m$,}
\\ \f m2V_kQ^{\f{k(m-2)}2}\cdot
n^rQ^{\f{km(n^r-1)}2}=\f{mn^r}2Q^{\f{k(mn^r-2)}2}V_k\mod{n^{r+1}A}&\t{if
$2\mid m$.}\endcases\endaligned$$
\par If $2\mid n$, by Lemma 2 we have
$$\f{U_{kmn^{s+1}}}{U_{kmn^s}}\e nQ^{\f{kmn^s(n-1)}2}\mod{|P^2-4Q|U_{\f{kmn^s}2}^2}\q (s\ge 1)$$
 and so
$U_{kmn^{s+1}}/U_{kmn^s}\e nQ^{\f{kmn^s(n-1)}2}\mod{n^2A}$ for $
s\ge 1.$ Thus
$$\f{U_{kmn^r}}{U_{kmn}}=\prod_{s=1}^{r-1}\f{U_{kmn^{s+1}}}{U_{kmn^s}}
\e
\prod_{s=1}^{r-1}\big(nQ^{\f{kmn^s(n-1)}2}\big)=n^{r-1}Q^{\f{km(n^r-n)}2}
\mod {n^rA}.$$ Note that $U_{2k}=U_kV_k$. By Lemma 2 we also have
$$\f{U_{kmn}}{U_k}\e mn\f{V_k}2Q^{k(\f {mn}2-1)}\mod{n^2A}.$$
Since $2\mid n$ and $n\mid U_k$ we see that $2\mid U_k$ and so
$2\mid V_k$ by (4). Hence
$$\f{U_{kmn^r}}{U_{k}}=\f{U_{kmn^r}}{U_{kmn}}\cdot
\f{U_{kmn}}{U_{k}}\e
\f{mV_k}2n^rQ^{\f{k(mn^r-2)}2}\mod{n^{r+1}A}.$$
 This completes the proof.

 \pro{Corollary 1} Let $k,m,n,r$ be positive integers, $P,Q\in\Bbb Z$,
 $(P,Q)=1$, $PQ(P^2-Q)(P^2-4Q)\not=0$ and $U_n=U_n(P,Q)$.
 If $p$ is an odd prime divisor of $U_k$, then
 $$\f{U_{kp^r}}{U_k}\e p^r\mod{p^{r+1}}.$$
\endpro
Proof. Taking $m=1$ and $n=p$ in Theorem 1 and then applying
Fermat's little theorem we get the result.
 \pro{Theorem 2} Let $k,m,n,r,s$ be positive integers, $P,Q\in\Bbb Z$,
 $(P,Q)=1$, $PQ(P^2-Q)(P^2-4Q)\not=0$ and $U_n=U_n(P,Q)$. Suppose
 $n^s\mid U_k$ and $n^{s+1}\nmid U_k$. Then
\par $\trm{i}$ $n^{r+s}\mid U_{kmn^r}$.
\par $\trm{ii}$ If $n$ is odd and
 $(m,n)=1$, then
  $n^{r+s+1}\nmid U_{kmn^r}$.\endpro
 Proof. Let $V_k=V_k(P,Q)$. If $2\mid n$, then $2\mid U_k$ and so $2\mid V_k$ by
 (4). Hence $n^r\mid\f{U_{kmn^r}}{U_k}$ by Theorem 1. That is
 $n^rU_k\mid U_{kmn^r}$. Since $n^s\mid U_k$ we must have
 $n^{r+s}\mid U_{kmn^r}$. This proves (i).
  \par Now consider (ii). Suppose $2\nmid n$, $(m,n)=1$ and $U_k=dn^s$. Then $n\nmid d$.
   From Theorem 1
we see that
$$U_{kmn^r}\e\cases
dmn^{r+s}Q^{\f{k(mn^r-1)}2}\mod{n^{r+s+1}}&\t{if $2\nmid
m$}, \\
\f{dm}2n^{r+s}Q^{\f{k(mn^r-2)}2}V_k\mod{n^{r+s+1}}&\t{if $2\mid
m$}.\endcases $$ Note that $n\nmid d,\ (m,n)=1$ and $(n,QV_k)=1$
by Lemma 4. By the above we obtain $n^{r+s+1}\nmid U_{kmn^r}$.
This finishes the proof.
\par Clearly Theorem 2 is a generalization of Lucas' law of
repetition.

\pro{Theorem 3} Let $P,Q\in\Bbb Z$, $PQ(P^2-Q)(P^2-4Q)\not=0$,
$(P,Q)=1$, $U_n=U_n(P,Q)$ and $V_n=V_n(P,Q)$.
\par $(\t{\rm i})$ If $p$ is an odd prime such that $p\mid U_m\
(m\ge 1)$, then
$$\ord_pU_m=\ord_pm+\ord_pU_{p-(\f{P^2-4Q}p)}.$$
\par $(\t{\rm ii})$ If $p$ is an odd prime such that $p\mid V_m\
(m\ge 1)$, then
$$\ord_pV_m=\ord_pm+\ord_pU_{p-(\f{P^2-4Q}p)}.$$
\endpro
Proof. If $p$ is an odd prime such that $p\mid U_m$ or $p\mid
V_m$, then $p\mid U_{2m}$. From Lemmas 4 and 5 we know that
$p\nmid Q$ and so $p\mid U_{p-(\f{P^2-4Q}p)}$.
 Let $r(p)$ be the least positive integer
$n$ such that $p\mid U_n$. Then $r(p)\mid p-(\f{P^2-4Q}p)$ and
$r(p)\mid 2m$ by Lemma 5. Hence $p\nmid r(p)$ and
$U_{p-(\f{P^2-4Q}p)}/U_{r(p)}\not\e 0\mod p$ by Theorem 1. Suppose
$2m=kp^tr(p)\ (p\nmid k)$. Then $t=\ord_p{2m}=\ord_pm$. Since
$V_n^2-(P^2-4Q)U_n^2=4Q^n$ by (4), we see that $(U_n,V_n)\mid
4Q^n$ and so $p\nmid (U_n,V_n)$ since $p\nmid Q$. Hence $p\mid
U_n$ implies $p\nmid V_n$, and $p\mid V_n$ implies $p\nmid U_n$.
Now applying Theorem 2 we get
$$\ord_p
U_{2m}=t+\ord_pU_{r(p)}=\ord_pm+\ord_pU_{p-(\f{P^2-4Q}p)}.$$ This
completes the proof. \pro{Corollary 2} Let $P\not=0$ be an
integer, $S_1=P^2+2$ and $S_{k+1}=S_k^2-2\ (k\ge 1)$. If $p$ is an
odd prime factor of $S_n\ (n\ge 1)$, then $p^{\alpha}\mid S_n$ if
and only if $p^{\alpha}\mid U_{p-(\f{P^2+4}p)}(P,-1)$.
\endpro
Proof. From (3) we see that $S_n=V_{2^n}(P,-1)$. Thus the result
follows from Theorem 3(ii). \par Let $F_n=U_n(1,-1)$ be the
Fibonacci sequence, and let $\{S_n\}$ be given by $S_1=3$ and
$S_{k+1}=S_k^2-2\ (k\ge 1)$. Recently R. McIntosh showed that
$p^2\nmid F_{p-(\f p5)}$ for any prime $p<10^{14}$. Thus it
follows from Corollary 2 that any square prime factor of $S_n$
should be greater than $10^{14}$. \pro{Corollary 3} Let
$a,n\in\Bbb Z$, $a\not=-1$ and $n\ge 1$. If $p$ is an odd prime
divisor of $a^n+1$, then
$$\ord_p(a^n+1)=\ord_pn+\ord_p(a^{p-1}-1).$$
\endpro
Proof. Since $p\mid a^n+1$ we must have $a\not=0,1$. By (2),
$V_n(a+1,a)=a^n+1$. So the result follows immediately from Theorem
3(ii).
\par From Corollary 3 we have
 \pro{Corollary 4} If $p$ is a prime divisor of
$2^{2^n}+1$, then $p^{\alpha}\mid 2^{2^n}+1$ if and only if
$p^{\alpha}\mid 2^{p-1}-1.$
\endpro
\par In the case $\alpha=2$, the result of Corollary 4 is well known.
 See [2,5,6].
 
\enddocument

\Refs
\ref\no 1\by L.E. Dickson\book History of the Theory of Numbers
\bookinfo Vol.I\publaddr Chelsea, New York\yr 1952\pages
393-407\endref \ref \no 2\by M. K$\check{\t{\rm
r}}$i$\check{\t{\rm z}}$ek, F. Luca and L. Somer\book 17 Lectures
on Fermat Numbers: From Number Theory to Geometry\publ
Springer-Verlag \publaddr New York\yr 2001\page 68\endref
 \ref\no 3\by P. Ribenboim\book The Book of
Prime Number Records\bookinfo 2nd ed.\publaddr Springer, Berlin\yr
1989\pages 44-50\endref \ref\no 4\by P. Ribenboim\book My numbers,
my friends\publ Springer-Verlag New York, Inc. \publaddr New York,
Berlin, London \yr 2000\pages 1-41\endref
 \ref\no 5\by P. Ribenboim\paper On the square factors of the
 numbers of Fermat and Ferentinou-Nicolacopoulou\jour Bull. Greek
 Math. Soc.\vol 20\yr 1979a\pages 81-92\finalinfo MR83f:10016.\endref
 \ref\no 6\by L.R.J. Warren and H.G. Bray\paper
On the square-freeness of Fermat and Mersenne numbers\jour Pacific
J. Math. \vol 22\yr 1967\pages 563-564\finalinfo
MR36$\#$3718\endref
\endRefs
\bye